\documentstyle{article}

\begin{document}

\newfont{\goth}{eufm10}

\setcounter{equation}{0}
\def\theequation{\thesection.\arabic{equation}}

\huge

\noindent
{\bf Intertwining Operators and}

\noindent
{\bf Quantum Homogeneous Spaces}
\vspace{5mm}

\normalsize

\noindent
{\bf Leonid L. Vaksman}
\footnote{Partially supported by the AMS fSU grant and by the ISF grant.}
\vspace{5mm}

\noindent
Physico-Technical Institute of Low Temperatures, Kharkiv, Ukraine
\vspace{7mm}

\noindent
{\bf Abstract.} In the present paper the algebras of functions
on quantum homogeneous spaces are studied. The author introduces
the algebras of kernels of intertwining integral operators and
constructs quantum analogues of the Poisson and Radon transforms
for some quantum homogeneous spaces. Some applications and the
relation to $q$-special functions are discussed.

\section{Introduction}

The integral transformations which intertwine quasi-regular
representations of groups are widely used in representation
theory, harmonic analysis, and mathematical physics. The
kernels of the intertwining integral operators are functions
on Cartesian products $X \times Y$ of $G$-spaces which are
constant on the $G$-orbits. These can be distributions, as
in the case of the Radon transform
\begin{equation}
\label{0.1}
K(x,y)=\delta(c-(x,y))
\end{equation}
and other operators related to the integral geometry.

The quantum group theory originated in the papers by
V.G.Drinfeld and M.Jimbo \cite{1,2} in the middle eighties
under the influence of the quantum inverse scattering method
developed by L.D.Faddeev and his school.

Our goal is to construct invariant kernels analogous to (\ref{0.1})
and the corresponding intertwining integral operators for a certain
class of quantum groups and their homogeneous spaces. At the same
time there will be defined real powers of the $q$-Poisson kernel
and shown their relation to the theory of $q$-special functions.
(Here $q \in (0,1)$ is a parameter related to the ``Plank constant''
$h$ used in \cite{1} by the equality $q=e^{-h/2}$.)

Quantum algebras of functions are deformations of algebras of
functions on Poisson manifolds. Given two Poisson manifolds $X$
and $Y$, the Poisson bracket on the Cartesian product $X \times Y$
can be defined by
\[ \{\varphi_{1} \otimes \psi_{1},\varphi_{2} \otimes \psi_{2}\}=
\{\varphi_{1},\varphi_{2}\} \otimes \psi_{1}\psi_{2}-
\varphi_{1}\varphi_{2} \otimes \{\psi_{1},\psi_{2}\} \]
where $\varphi_{i} \in C^{\infty}(X)$, $\psi_{i} \in C^{\infty}(Y)$
($i=1,2$). The 'minus' sign is important when $X$ and $Y$ are
Poisson homogeneous $G$-spaces for a Poisson Lie group $G$
\cite{1}. Due to it, the Poisson bracket of $G$-invariants
in $C^{\infty}(X \times Y)$ is again a $G$-invariant. The quantum
analogue of this simple observation is that the intertwining
kernels form an algebra (see Section 2). In particular, a polynomial
of an intertwining kernel is again an intertwining kernel. However,
as we see on the example of (\ref{0.1}), it does not suffice to have
polynomials only, we need distributions. The corresponding constructions
are given in Sections 4 and 5.

Finally, it is important to note that the algebras of intertwining
kernels are described by extremely simple commutation relations
betwee the generators (see Section 6). This can be explained by
the fact that quantization preserves the dimensions of the graded
components of the space of intertwining operators.

\section{Integral Operators and $A$-module algebras}

\setcounter{equation}{0}
\def\theequation{\thesection.\arabic{equation}}

Consider a Hopf algebra $A$ over {\bf C}. Recall (cf. \cite{1}) that
the comultiplication $\Delta:A \rightarrow A \otimes A$ is a
homomorphism, which allows to introduce the operation of tensor
product of $A$-modules by. Also, given an $A$-module $V$, one can
define the dual $A$-module $V^{*}$ by
\[ (av^{*})(v)=v^{*}(S(a)v) \]
where $a \in A$, $v \in V$, $v^{*} \in V^{*}$, and $S:A \rightarrow A$
is the antipode. The trivial $A$-module {\bf C} is defined by the
counit $\varepsilon:A \rightarrow {\bf C}$. By definition, the
$A$-invariants of an $A$-module $V$ are the vectors $v \in V$ such
that $av=\varepsilon(a)v$ for any $a \in A$. Recall, finally, that
the antipode $S$ is both an algebra and a coalgebra anti-isomorphism
of $A$.

Important examples of Hopf algebra are the quantum universal enveloping
algebras $U_{q}\mbox{\goth g}$ described in \cite{1,2}. Given a simple
complex Lie algebra {\goth g} with Cartan matrix $A=(a_{ij})_{i,j=1}^{n}$
where $n={\rm rank} \, {\goth g}$, the algebra structure on
$U_{q}\mbox{\goth g}$ is defined by the generators $X_{i}^{\mp}$,
$K_{i}^{\pm}$ ($i=1,...,n$) and the relations
\begin{eqnarray}
\label{1.1-1}
K_{i}^{\pm}K_{j}^{\pm}=K_{j}^{\pm}K_{i}^{\pm}, &
K_{i}^{+}K_{i}^{-}=K_{i}^{-}K_{i}^{+}=1, \\
\label{1.1-2}
K_{i}^{+}X_{j}^{\pm}=q^{\pm \frac{a_{ij}}{2}}X_{j}^{\pm}K_{i}^{+}, &
K_{i}^{-}X_{j}^{\pm}=q^{\mp \frac{a_{ij}}{2}}X_{j}^{\pm}K_{i}^{-},
\end{eqnarray}

\vspace{-5.6mm}

\begin{eqnarray}
\label{1.1-3}
X_{i}^{+}X_{j}^{-}-X_{j}^{-}X_{i}^{+}=\delta_{ij}
\frac{K_{i}^{+2}-K_{i}^{-2}}{q-q^{-1}}, \hspace{19mm} \\
\label{1.1-4}
\sum_{k=0}^{1-a_{ij}} (-1)^{k} q_{i}^{k(1-a_{ij}-k)} \left[
\begin{array}{c}
1-a_{ij} \\ k
\end{array}
\right] _{q_{i}^{2}} (X^{\pm}_{i})^{1-a_{ij}-k}
X^{\pm}_{j}(X^{\pm}_{i})^{k}=0 \\
\nonumber
\mbox{\rm when} \; \; i \not = j \hspace{40mm}
\end{eqnarray}
where $q_{i}=q^{\frac{(\alpha_{i},\alpha_{i})}{2}}$ (with
$\{\alpha_{i}\}_{i=1}^{n}$ being a system of simple coroots).
Here $\left[
\begin{array}{c}
m \\ n
\end{array}
\right]_{t}=\frac{(t;t)_{m}}{(t;t)_{n}(t;t)_{m-n}}$ is the
$t$-binomial coefficient, $(a;t)_{k}=(1-a)(1-at)...(1-at^{k-1})$.

The coalgebra structure on $U_{q}\mbox{\goth g}$ is given by
\begin{equation}
\label{1.1-5}
\Delta(K_{i}^{\pm})=K_{i}^{\pm} \otimes K_{i}^{\pm},
\; \; \; \Delta(X_{i}^{\pm})=K_{i}^{-} \otimes
X_{i}^{\pm}+X_{i}^{\pm} \otimes K_{i}^{+} \\
\end{equation}
\begin{equation}
\label{1.1-6}
\varepsilon(K_{i}^{\pm})=1, \; \; \varepsilon(X_{i}^{\pm})=0, \; \;
S(K_{i}^{\pm})=K_{i}^{\mp}, \; \; S(X_{i}^{\pm})=-q_{i}^{\mp 1}X_{i}^{\pm}
\end{equation}

\noindent
{\bf Definition 2.1.} Let $A$ be a Hopf algebra and $F$ an
algebra equipped with an $A$-module structure. We say that
$F$ is an $A$-{\em module algebra\/} if the multiplication
map
\[ m:F \otimes F \rightarrow F, \; \; \;
m(f_{1} \otimes f_{2})=f_{1}f_{2} \]
is a morphism of $A$-modules.
\vspace{4mm}

As follows from (\ref{1.1-5})-(\ref{1.1-6}), in the case of
$A=U_{q}\mbox{\goth g}$ the condition of Definition 2.1
means that
\begin{equation}
\label{1.2}
K_{i}^{\pm}(f_{1}f_{2})=K_{i}^{\pm}(f_{1})K_{i}^{\pm}(f_{2})
\end{equation}
\begin{equation}
\label{1.3}
X_{i}^{\pm}(f_{1}f_{2})=K_{i}^{-}(f_{1})X_{i}^{\pm}(f_{2})
+X_{i}^{\pm}(f_{1})K_{i}^{+}(f_{2})
\end{equation}
for any $i=1,...,n$, $f_{1},f_{2} \in F$. The last equality
is an analogue of the Leibnitz formula.

Reasonable examples of $A$-module algebras -- algebras of functions
on quantum homogeneous spaces -- are given in the paper \cite{3} by
L.D.Faddeev, N.Yu.Re\-shetikhin, L.A.Takhtajan. One more example is
the tensor algebra $T$ generated by $A$-modules $V_{1}, V_{1}^{*},
V_{2}, V_{2}^{*}, ..., V_{k}, V_{k}^{*}$. The algebras in \cite{3}
are quotients of the tensor algebras over ideals defined in terms
of certain solutions of the quantum Yang-Baxter equation.

As in the limit case $q=1$, $A$-module morphisms can be obtained
by using the operation of tensor coevaluation (for example, $V_{j}^{*}
\otimes V_{j} \rightarrow {\bf C}$), and invariants can be obtained by
using the canonical embeddings ${\bf C} \hookrightarrow V_{j} \otimes
V_{j}^{*}$. We will need the following well-known properties of
$A$-modules.
\vspace{4mm}

\noindent
{\bf Proposition 2.2.} {\em Let\/} $U$ {\em and\/} $V$ {\em be\/}
$A$-{\em modules, and\/} $l:V \otimes U \rightarrow {\bf C}$
{\em a morphism of\/} $A$-{\em modules. Then\/} $l(av \otimes u)=
l(v \otimes S(a)u)$ {\em for any\/} $a \in A$, $u \in U$, $v \in V$.
\vspace{4mm}

\noindent
{\em Proof\/}. In the special case when $V=U^{*}$ the statement
of the proposition follows from the definition of dual module.
The general case can be reduced to this special case, since the
linear operator $V \rightarrow U^{*}$ defined by $l$ is a morphism
of $A$-modules:
\[ V \simeq V \otimes {\bf C} \hookrightarrow V \otimes U
\otimes U^{*} \stackrel{l \otimes {\rm id}}{\longrightarrow}
{\bf C} \otimes U^{*} \simeq U^{*} \]

\noindent
{\bf Corollary 2.3.} {\em Let\/} $F$ {\em be an\/} $A$-{\em module
algebra and\/} $\nu:F \rightarrow {\bf C}$ {\em a morphism of\/}
$A$-{\em modules\/} ({\em called\/} invariant integral {\em on\/} $F$:
$\int fd\nu \stackrel{\rm def}{=} \nu(f)$). {\em Then for any\/}
$a \in A$ {\em and\/} $f_{1},f_{2} \in F$, {\em the following
``integration by parts'' formula holds\/}:
\begin{equation}
\label{1.4}
\int (af_{1})f_{2}d\nu=\int f_{1}(S(a)f_{2})d\nu
\end{equation}

\noindent
{\bf Proposition 2.4.} {\em Consider\/} $A$-{\em modules\/} $V_{1},
V_{2}$ {\em and the image of\/} $Hom_{A}(V_{1},V_{2})$ {\em under
the canonical isomorphism\/} $Hom_{\bf C}(V_{1},V_{2}) \simeq V_{2}
\otimes V_{1}^{*}$.

(i) {\em Every invariant in\/} $V_{2} \otimes V_{1}^{*}$ {\em belongs
to the image of\/} $Hom_{A}(V_{1},V_{2})$.

(ii) {\em Given an element\/} $u$ {\em of the image of\/}
$Hom_{A}(V_{1},V_{2})$, {\em one has that\/}
\begin{equation}
\label{1.5}
(a \otimes 1)u=(1 \otimes S^{-1}(a))u
\end{equation}
{\em for any\/} $a \in A$.

(iii) {\em If\/} $u$ {\em satisfies\/} (\ref{1.5}) {\em for any\/}
$a \in A$, {\em then\/} $u$ {\em is an invariant\/}.
\vspace{4mm}

\noindent
{\em Proof\/}. If $u \in V_{2} \otimes V_{1}^{*}$ is an invariant,
then the multiplication by it $V_{1} \rightarrow V_{2} \otimes
V_{1}^{*} \otimes V_{1}$ with the subsequent coevaluation $V_{2}
\otimes V_{1}^{*} \otimes V_{1} \rightarrow V_{2} \otimes {\bf C}
\simeq V_{2}$ gives an element of $Hom_{A}(V_{1},V_{2})$.
This proves the first statement, while the other two follow from
the definitions.
\vspace{4mm}

\noindent
{\bf Corollary 2.5.} {\em The space of invariants of\/} $V_{2}
\otimes V_{1}^{*}$ {\em is determined by the condition\/}
(\ref{1.5}) {\em and canonically isomorphic to\/} $Hom_{A}
(V_{1},V_{2})$.
\vspace{4mm}

Now we are going to give a definition of the algebra of
intertwining kernels. To any algebra $F$ we assign an algebra
$F^{\rm op}$ obtained from $F$ by changing the multiplication
by the opposite one: $m^{\rm op}(f_{1} \otimes f_{2})=f_{2}f_{1}$.
We have the following proposition.
\vspace{4mm}

\noindent
{\bf Proposition 2.6.} {\em For any\/} $A$-{\em module algebras\/}
$F_{1}$ {\em and\/} $F_{2}$, {\em the invariants in\/} $F_{2}
\otimes F_{1}^{\rm op}$ {\em form a subalgebra\/}.
\vspace{4mm}

\noindent
{\em Proof\/}. By Corollary 2.5, we get a system of equations on
the invariants:
\begin{equation}
\label{1.6}
\forall a \in A: \; \; (a \otimes 1)f=1 \otimes S^{-1}(a)f
\end{equation}

Equip the algebra $F_{2} \otimes F_{1}^{\rm op}$ with an $A \otimes
A^{\rm op}$-module algebra structure, letting
\[ a_{2} \otimes a_{1}:f_{2} \otimes f_{1} \mapsto
a_{2}f_{2} \otimes S^{-1}(a_{1})f_{1} \]
for any $a_{1},a_{2} \in A$ and $f_{1} \in F_{1}^{\rm op}$, $f_{2}
\in F_{2}$. The system (\ref{1.6}) takes the form
\[ (a \otimes 1)f=(1 \otimes a)f \]

It is clear that its solutions form a subalgebra.
\vspace{4mm}

\noindent
{\bf Definition 2.7.} The subalgebra of invariants in $F_{2} \otimes
F_{1}^{\rm op}$ is called {\em algebra of intertwining kernels\/}.
(While the whole algebra $F_{2} \otimes F_{1}^{\rm op}$ is called
{\em algebra of kernels\/}.)
\vspace{4mm}

If there exists an invariant integral $\nu:F_{1} \rightarrow {\bf C}$,
then the map $i:F_{1} \rightarrow F_{1}^{*}$ given by
\[ (if_{1})(f_{2})=\int f_{1}f_{2}d\nu \]
is a morphism of $A$-modules. This allows to assign a morphism
$F_{1} \rightarrow F_{2}$ to any intertwining kernel $K$. In
other words, intertwining kernels correspond to the integral
operators
\[ f_{1} \mapsto ({\rm id} \otimes \nu)(K(1 \otimes f_{1})) \]
which intertwine the representations of $A$ in $F_{1}$ and $F_{2}$
(the functions under the symbol of integral are being multiplied in
the algebra $F_{2} \otimes F_{1}$, not in $F_{2} \otimes F_{1}^{\rm op}$).

Note that the replacement of $F_{1}$ by $F_{1}^{\rm op}$ results in
the change of sign in the commutators, and in the limit as $q \rightarrow
1$ it yeilds the above-mentioned change of sign in the formula for the
Poisson bracket on the Cartesian product of the corresponding Poisson
homogeneous spaces.
\vspace{4mm}

Consider the case $\mbox{\goth g}=\mbox{\goth sl}(n+1)$. The algebra
$U_{q}\mbox{\goth sl}(n+1)$ can be equipped with a anti-linear
anti-involution given by
\[ \left(K_{i}^{\pm}\right)^{*}=K_{i}^{\mp}, \; \; \;
\left(X_{i}^{\pm}\right)^{*}=(-1)^{\delta_{i,1}}X_{i}^{\mp} \]

This is an algebra anti-automorphism and a coalgebra automorphism.
Moreover, the map $a \mapsto (S(a))^{*}$ is an involution. That is,
the pair $(U_{q}\mbox{\goth sl}(n+1),*)$ denoted by $U_{q}
\mbox{\goth su}(n,1)$ is a Hopf $*$-algebra. The notation
$U_{q}\mbox{\goth su}(n,1)$ is justified by the fact that in
the limit case $q=1$ the $*$-representations of this algebra
are related to the unitary representations of $SU(n,1)$.
\vspace{4mm}

\noindent
{\bf Definition 2.8.} Let $A$ be a Hopf $*$-algebra. Suppose that
a $*$-algebra $F$ is equipped with an $A$-module structure so that
the multiplication map $m:F \otimes F \rightarrow F$ is a morphism
of $A$-modules (that is, $F$ is an $A$-module algebra).

We say that $F$ is an $A$-{\em module\/} $*$-{\em algebra\/} if the
following condition is satisfied:
\begin{equation}
\label{1.7}
(af)^{*}=(S(a))^{*}f^{*}
\end{equation}
for any $a \in A$ and $f \in F$.
\vspace{4mm}

Given an $A$-module $*$-algebra $F$, we will choose an
invariant integral $\nu:F \rightarrow {\bf C}$ so that
\[ \int f^{*}d\nu=\overline{\int fd\nu}, \; \; \;
\int f^{*}fd\nu>0 \; \; {\rm if} \; \; f \not =0 \]

Then we introduce a scalar product in $F$ by letting
\[ (f_{1},f_{2})=\int f_{2}^{*}f_{1}d\nu \]

It is easy to see that (\ref{1.7}) is equivalent to the condition
that the quasi-regular representation of $A$ in $F$ given by
${\cal R}(a)f=af$ is a $*$-representation. Indeed, $(af_{1},f_{2})
=\int \left(S^{-1}(a)f_{2}^{*}\right)f_{1}d\nu=\int(a^{*}f_{2})^{*}
f_{1}d\nu=(f_{1},a^{*}f_{2})$.

Suppose that we have two $A$-module $*$-algebras $F_{1}$ and $F_{2}$.
It is easy to show that there exists a unique anti-linear anti-involution
in the algebra of kernels such that
\[ \int K^{*}(1 \otimes f)d\nu=\left( \int K \left(1 \otimes
f^{*} \right) d\nu \right)^{*} \]
for any $a \in A$, $f \in F_{1}$, and $K \in F_{2} \otimes F_{1}^{*}$.
In other words, real kernels $K=K^{*}$ correspond to real integral
operators.

If $K$ is an intertwining kernel, then for any $a \in A$
and $f \in F_{1}$, we get that
\begin{eqnarray*}
\int K^{*}(1 \otimes af)d\nu & = & \left( \int K \left(1 \otimes
(af)^{*} \right) d\nu \right)^{*}= \\
=\left( \int K \left(1 \otimes (S(a))^{*}f^{*} \right) d\nu \right)^{*}
 & = & \left( (S(a))^{*} \int K \left( 1 \otimes f^{*} \right)
d\nu \right)^{*}= \\
=a \left( \int K \left( 1 \otimes f^{*} \right)
d\nu \right)^{*} & = & a \int K^{*}(1 \otimes f)d\nu
\end{eqnarray*}
Therefore, $K^{*}$ is also an intertwining kernel.

Given a kernel $K$, one can find $K^{*}$, using explicit formulas for
the invariant integral. If the algebra $F_{1}$ has an exact irreducible
$*$-representation $\pi$, usually there is a non-negative element
$Q \in F_{1}$ such that
\[ \int fd\nu=tr \, \pi(fQ) \]
for any $f \in F_{1}$ (see (\ref{2.8})). In the case of a reducible
representation $\pi$ the usual trace is replaced by the trace on
the corresponding von Neumann algebra. We have that
\begin{equation}
\label{1.8}
K^{*}=(1 \otimes Q)^{-1}K^{* \otimes *}(1 \otimes Q)
\end{equation}
In all the examples considered below in the paper the elements
$Q$ and $Q^{-1}$ does not belong in fact to the algebra $F_{1}$
of basic functions but to a space of distributions (see Section 3).

Following \cite{3}, we describe a $U_{q}\mbox{\goth su}(n,1)$-module
$*$-algebra ${\bf C}_{q,q^{-1}}^{n+1}$ which is a quantum analogue
of the algebra of polynomial functions on a vector space \footnote{In
\cite{3} this algebra was denoted by ${\bf C}_{q,q^{-1}}^{n+1}$.}.
As a $*$-algebra with unit element, ${\bf C}_{q,q^{-1}}^{n+1}$ is
generated by $\{z_{j}\}_{j=0}^{n}$ and the relations
\begin{equation}
\label{1.9}
z_{i}z_{j}=qz_{j}z_{i} \; \; (i<j), \; \; \;
z_{i}z_{j}^{*}=qz_{j}^{*}z_{i} \; \; (i \not =j)
\end{equation}
\begin{equation}
\label{1.10}
z_{i}z_{i}^{*}-z_{i}^{*}z_{i}=(-1)^{\delta_{i,0}}
\left(q^{-2}-1\right)\sum_{k>i}z_{k}z_{k}^{*}
\end{equation}

Define an action of $U_{q}\mbox{\goth su}(n,1)$ on the generators
$z_{j}$ by
\begin{eqnarray}
\label{1.12}
X_{i}^{+}z_{j}=\delta_{ij}z_{j-1} & (j \not =0),
& X_{i}^{+}z_{0}=0 \\
\label{1.13}
X_{i}^{-}z_{j}=\delta_{i,j+1}z_{j+1}  &
(j \not =n), & X_{i}^{-}z_{n}=0
\end{eqnarray}
\begin{equation}
\label{1.14}
K_{i}^{\pm}z_{j}=\left( \delta_{i,j+1}q^{\pm \frac{1}{2}}+
\delta_{ij}q^{\mp \frac{1}{2}} \right)z_{j}
\end{equation}

This acton can be extended to an action on the whole
${\bf C}_{q,q^{-1}}^{n+1}$ by (\ref{1.2}), (\ref{1.3},
and (\ref{1.7}). For instance, we get that
\begin{eqnarray}
\label{1.15}
X_{i}^{+}z_{j}^{*}=-q^{-1}\delta_{i,j+1}z_{j+1}^{*} & (0<j<n) \\
\label{1.16}
X_{i}^{-}z_{j}^{*}=q\delta_{ij}z_{j-1}^{*} & (j>1)
\end{eqnarray}
\begin{equation}
\label{1.17}
K_{i}^{\pm}z_{j}^{*}=\left( \delta_{i,j+1}q^{\mp
\frac{1}{2}}+\delta_{ij}q^{\pm \frac{1}{2}} \right)z_{j}^{*}
\end{equation}

We will denote the generators of the algebra ${\bf C}_{q,q^{-1}}^{n+1}
\otimes {\bf C}_{q,q^{-1}}^{n+1}$ by $z_{j}=z_{j} \otimes 1$ and
$\zeta_{j}=1 \otimes z_{j}$. As vector spaces, the algebras
${\bf C}_{q,q^{-1}}^{n+1} \otimes {\bf C}_{q,q^{-1}}^{n+1}$ and
${\bf C}_{q,q^{-1}}^{n+1} \otimes \left( {\bf C}_{q,q^{-1}}^{n+1}
\right)^{\rm op}$ coincide.

Consider the kernels $t,\tau,K_{1},K_{2} \in {\bf C}_{q,q^{-1}}^{n+1}
\otimes \left( {\bf C}_{q,q^{-1}}^{n+1} \right)^{\rm op}$ given by
\[ t=z_{0}z_{0}^{*}-\sum_{j=1}^{n}z_{j}z_{j}^{*}, \; \; \;
\tau=\zeta_{0}\zeta_{0}^{*}-\sum_{j=1}^{n}\zeta_{j}\zeta_{j}^{*}, \]
\[ K_{1}=z_{0}\zeta_{0}^{*}-\sum^{n}_{j=1}z_{j}\zeta_{j}^{*}, \; \; \;
K_{2}=z_{0}^{*}\zeta_{0}-\sum^{n}_{j=1}z_{j}^{*}\zeta_{j} \]

\noindent
{\bf Proposition 2.9.} {\em The kernels\/} $t,\tau,K_{1},K_{2}$
{\em are intertwining. Moreover\/}, $t$ {\em and\/} $\tau$
{\em belong to the center of\/} ${\bf C}_{q,q^{-1}}^{n+1}
\otimes \left( {\bf C}_{q,q^{-1}}^{n+1} \right)^{\rm op}$
{\em and\/}
\begin{equation}
\label{1.18}
K_{1}K_{2}-q^{2}K_{2}K_{1}-\left(1-q^{2}\right)t\tau=0
\end{equation}

\noindent
{\bf Corollary 2.10.} $(K_{1}K_{2})^{m}$ {\em is an intertwining
kernel for any\/} $m \in {\bf Z}_{+}$.
\vspace{4mm}

In the theory of representations of $SU(n,1)$ an important role
is played by non-integer powers of the Poisson kernel $P=K_{1}K_{2}$.
Their quantum analogues will be defined in Section 4.

\section{Distributions and Intertwining Kernels}

\setcounter{equation}{0}
\def\theequation{\thesection.\arabic{equation}}

Recall (cf. \cite{4}) the definitions of the algebras of regular
functions on a quantum cone ${\bf C}\left[C^{2n+1}\right]_{q}$
and a quantum hyperboloid ${\bf C}\left[H^{2n+1}_{+}\right]_{q}$.
These are the quotients of ${\bf C}^{n+1}_{q,q^{-1}}$ over the
ideals generated by the central elements
\[ f_{0}=z_{0}z_{0}^{*}-\sum_{k=1}^{n}z_{k}z_{k}^{*}
\; \; \; \mbox{\rm and} \; \; \; f_{1}=f_{0}-1 \; \; \;
\mbox{\rm respectively} \]

Since the elements $f_{0}$ and $f_{1}$ are invariants and
$f_{0}=f_{0}^{*}$, $f_{1}=f_{1}^{*}$, the algebras ${\bf C}
\left[C^{2n+1}\right]_{q}$ and ${\bf C}\left[H^{2n+1}_{+}\right]_{q}$
inherit a $U_{q}\mbox{\goth su}(n,1)$-module $*$-algebra structure.

However, the algebra ${\bf C}\left[H^{2n+1}_{+}\right]_{q} \otimes
{\bf C}\left[H^{2n+1}_{+}\right]_{q}^{\rm op}$ contains only polynomial
intertwining kernels. In order to define a space of distributions
we first consider some corollaries of the commutation relations
(\ref{1.9}),(\ref{1.10}).
\vspace{4mm}

\noindent
{\bf Proposition 3.1.} {\em Let\/} $x_{i} \in {\bf C}_{q,q^{-1}}^{n+1}$
{\em be given by\/}
\begin{equation}
\label{2.1}
x_{0}=z_{0}z_{0}^{*}-\sum_{k>0}z_{k}z_{k}^{*},
\; \; \; x_{j}=\sum_{k \geq j}z_{k}z_{k}^{*}
\; \; (j>0)
\end{equation}
{\em Then we have that\/}
\begin{eqnarray*}
x_{i}x_{j} & = & x_{j}x_{i} \\
z_{i}x_{j} & = & x_{j}z_{i} \hspace{5mm} \; (i \geq j) \\
z_{i}x_{j} & = & q^{2}x_{j}z_{i} \; \; \: (i<j)
\end{eqnarray*}

\noindent
{\bf Corollary 3.2.} {\em If\/} $f$ {\em is an element of either\/}
${\bf C}\left[C^{2n+1}\right]_{q}$ {\em or\/} ${\bf C}\left[H^{2n+1}_{+}
\right]_{q}$, {\em then there exists a unique decomposition\/}
\begin{equation}
\label{2.2}
f=\sum_{I \cdot J=0} \left(z^{*}\right)^{I}
f_{IJ}(x_{1},...,x_{n})z^{J}
\end{equation}
{\em where\/} $I=(i_{0},...,i_{n})$, $J=(j_{0},...,j_{n})$,
$I \cdot J=(i_{0}j_{0},...,i_{n}j_{n})$, and $\left(z^{*}\right)^{I}
=\left(z_{0}^{*}\right)^{i_{0}}\left(z_{1}^{*}\right)^{i_{1}}...
\left(z_{n}^{*}\right)^{i_{n}}$, $z^{J}=z_{0}^{j_{0}}z_{1}^{j_{1}}
...z_{n}^{j_{n}}$.
\vspace{4mm}

Of course, here we have preserved the same notation for the
generators of the quotient algebras, the sum in (\ref{2.2})
is finite, and the coefficients $f_{IJ}$ are polynomials.

The action of the generators $X_{i}^{\pm}$ and $K_{i}^{\pm}$
of $U_{q}\mbox{\goth su}(n,1)$ on the terms of the sum in
(\ref{2.2}) is described by $q$-difference operators, as
follows from (\ref{1.2}), (\ref{1.3}), and the following
simple statement.
\vspace{4mm}

\noindent
{\bf Proposition 3.3.} {\em For any polynomial\/} $f(x_{1},...,x_{n})$
{\em in the algebra\/} ${\bf C}_{q,q^{-1}}^{n+1}$, {\em we have that\/}
\begin{equation}
\label{2.3}
X_{i}^{+}f=q^{\frac{3}{2}}z_{i}^{*}(B_{i}f)(x_{1},...,x_{i},
q^{2}x_{i+1},...,q^{2}x_{n})z_{i-1}
\end{equation}
\begin{equation}
\label{2.4}
X_{i}^{-}f=-(-1)^{\delta_{i,1}}q^{\frac{3}{2}}z_{i-1}^{*}
(B_{i}f)(x_{1},...,x_{i},q^{2}x_{i+1},...,q^{2}x_{n})z_{i}
\end{equation}
\begin{equation}
\label{2.6}
K_{i}^{\pm}f=f
\end{equation}
{\em where\/}
\[ B:f(t) \mapsto \frac{f\left(q^{2}t\right)-f(t)}{q^{2}t-t}, \]
\[ B_{i}=\underbrace{\mbox{\rm id} \otimes ... \otimes
\mbox{\rm id}}_{i-1} \otimes B \otimes \mbox{\rm id}
\otimes ... \otimes \mbox{\rm id} \]

The equalities (\ref{2.3})-(\ref{2.6}) and the commutation relations
\begin{eqnarray*}
z_{i}f(x_{1},...,x_{n}) & = & f(x_{1},...,x_{i},
q^{2}x_{i+1},...,q^{2}x_{n})z_{i} \\
z_{i}^{*}f(x_{1},...,x_{n}) & = & f(x_{1},...,x_{i},
q^{-2}x_{i+1},...,q^{-2}x_{n})z_{i}^{*}
\end{eqnarray*}
enable us to extend the ring of polynomials ${\bf C}[x_{1},...,x_{n}]$
whose elements have been so far the coefficients $f_{IJ}$ of the
sum in (\ref{2.2}). Namely, let
\[ \mbox{\goth M}_{\beta} \stackrel{\rm def}{=} \]
\[ \left\{\left(q^{2(m_{1}+\beta)},...,
q^{2(m_{n}+\beta)}\right) \in {\bf R}^{n} \;
| \; (m_{1},...,m_{n}) \in {\bf Z}^{n}, \:
m_{1}<m_{2}<...<m_{n} \right\} \]
where $0 \leq \beta<1$.

Consider the algebra ${\cal D}_{\beta}$ of functions on ${\bf R}^{n}$
with finite support ($\#(\mbox{\rm supp} \, f)<\infty$) such that
$\mbox{\rm supp} \, f \subset \mbox{\goth M}_{\beta}$. We call {\em
basic functions\/} on the quantum hyperboloid sums of the form
(\ref{2.2}) where $f_{IJ} \in {\cal D}_{\beta}$. Similarly we define
basic functions on the quantum cone. The spaces of basic functions
will be denoted by ${\cal D}\left(H^{2n+1}_{+}\right)_{q,\beta}$
and ${\cal D}\left(C^{2n+1}\right)_{q,\beta}$ respectively. It is
clear that they can be equipped with a $U_{q}\mbox{\goth su}(n,1)$-module
$*$-algebra structure by (\ref{1.12})-(\ref{1.17}) and
(\ref{2.3})-(\ref{2.6}).

Let $\left\{f^{(m)}\right\}$ be a sequence of elements of ${\bf C}
\left[H^{2n+1}_{+}\right]_{q}$ and $f \in {\cal D}\left(H^{2n+1}_{+}
\right)_{q,\beta}$. We will say that the sequence $f^{(m)}$ converges
to $f$ if for any $I,J$ we have pointwise on $\mbox{\goth M}_{\beta}$
that
\[ \lim_{m \rightarrow \infty} f_{IJ}^{(m)}=f_{IJ} \; \; \;
\mbox{\rm and} \; \; \; \#\left(\bigcup_{m}\left\{(I,J) \; | \;
f_{IJ}^{(m)} \not \equiv 0 \right\}\right) < \infty \]

The same way we define convergence of a sequence of elements of
${\bf C}\left[C^{2n+1}\right]_{q}$ to an element of ${\cal D}
\left(C^{2n+1}\right)_{q,\beta}$.
\vspace{4mm}

\noindent
{\bf Proposition 3.4.} (i) {\em The multiplication and the action of\/}
$U_{q}\mbox{\goth su}(n,1)$ {\em can be extended continuously from\/}
${\bf C}\left[C^{2n+1}\right]_{q}$ {\em to\/} ${\cal D}\left(C^{2n+1}
\right)_{q,\beta}$ {\em and from\/} ${\bf C}\left[H^{2n+1}_{+}\right]_{q}$
{\em to\/} ${\cal D}\left(H^{2n+1}_{+}\right)_{q,\beta}$.

(ii) {\em The formula\/}
\begin{equation}
\label{2.7}
\int fd\nu_{\beta}=\left(1-q^{2}\right)^{n} \sum_{(x_{1},...,x_{n})
\in \mbox{\goth M}_{\beta}} f_{00}(x_{1},...,x_{n})x_{1}...x_{n}
\end{equation}
{\em defines an invariant integral on both\/} ${\cal D}\left(C^{2n+1}
\right)_{q,\beta}$ {\em and\/} ${\cal D}\left(H^{2n+1}_{+}\right)_{q,\beta}$.

(iii) {\em The involution\/} $*$ {\em can be extended
continuously to\/} ${\cal D}\left(C^{2n+1}\right)_{q,\beta}$ {\em and\/}
${\cal D}\left(H^{2n+1}_{+}\right)_{q,\beta}$. {\em Then we have that\/}
$\int f^{*}fd\nu_{\beta}>0$ {\em for any\/} $f \not =0$.
\vspace{4mm}

\noindent
{\em Proof\/}. When checking the statements (i) and (ii) the Hopf algebra
$U_{q}\mbox{\goth sl}(n+1)$ can be replaced by a Hopf subalgebra
isomorphic to $U_{q}\mbox{\goth sl}(2)$ and generated by the elements
$X_{i}^{\pm}$ and $K_{i}^{\pm}$, with $i$ being fixed. Then it suffices
to use (\ref{1.12})-(\ref{1.17}) and (\ref{2.3})-(\ref{2.6}).

Let us prove (iii). The involution can be extended in an obvious way.
The positivity of the integral (\ref{2.7}) follows from the fact that
\begin{equation}
\label{2.8}
\int f d\nu=\left(1-q^{2}\right)^{n} \frac{1}{2\pi} \int
tr \, \pi_{\varepsilon}^{(\beta,\varphi)}(fx_{1}...x_{n})
\end{equation}
where in the case $\varepsilon=0$ $\pi_{0}^{(\beta,\varphi)}$ is a
certain $*$-representation of ${\cal D}\left(C^{2n+1}\right)_{q,\beta}$
and in the case $\varepsilon=1$ $\pi_{1}^{(\beta,\varphi)}$ is a
certain $*$-representation of ${\cal D}\left(H^{2n+1}_{+}\right)_{q,\beta}$.

The representation $\pi_{\varepsilon}^{(\beta,\varphi)}$ is defined in
the space of functions $\psi$ with finite support on ${\bf Z} \times
{\bf Z}_{+}^{n-1}$. The scalar product in this space is given in a
standard way by
\[ (\psi_{1},\psi_{2})=\sum_{(j_{1},...,j_{n})}
\psi_{1}(j_{1},...,j_{n}) \overline{\psi_{2}(j_{1},...,j_{n})} \]
The operators $\pi_{\varepsilon}^{(\beta,\varphi)}(z_{k})$ and
$\pi_{\varepsilon}^{(\beta,\varphi)}(f(x_{1},...,x_{n}))$ are
given by
\begin{eqnarray}
\label{pi-1}
\left(\pi_{\varepsilon}^{(\beta,\varphi)}(z_{n})\psi\right)
(j_{1},...,j_{n}) & = &
e^{-i\varphi}q^{\beta+j_{1}+...+j_{n}}\psi(j_{1},...,j_{n}) \\
\label{pi-2}
\left(\pi_{\varepsilon}^{(\beta,\varphi)}(z_{k})\psi\right)
(j_{1},...,j_{n}) & = &
e^{-i\varphi}q^{\beta+j_{1}+...+j_{k}}\left(1-q^{2j_{k+1}}\right) \cdot \\
\nonumber
   & \cdot & \psi(j_{1},...,j_{k},j_{k+1}+1,j_{k+2},...,j_{n}), \; \; \:
(0<k<n) \\
\label{pi-3}
\left(\pi_{\varepsilon}^{(\beta,\varphi)}(z_{0})\psi\right)
(j_{1},...,j_{n}) & = &
q^{\beta}\left(\varepsilon+q^{2(j_{1}+\beta)}\right)\psi(j_{1}+1,
j_{2},...,j_{n}) \\
\label{pi-4}
\left(\pi_{\varepsilon}^{(\beta,\varphi)}(f)\psi\right)
(j_{1},...,j_{n}) & = & f\left(q^{2(\beta+j_{1})},...,
q^{2(\beta+j_{1}+...+j_{n})}\right)\psi(j_{1},...,j_{n})
\end{eqnarray}

Then the positivity of the integral follows from the exactness of the
representations $\oplus \int \pi_{0}^{(\beta,\varphi)}d\varphi$,
$\oplus \int \pi_{1}^{(\beta,\varphi)}d\varphi$.
\vspace{4mm}

\noindent
{\em Remarks\/}. (i) In what follows the algebras ${\cal D}\left(
C^{2n+1}\right)_{q,\beta}$ and ${\cal D}\left(H^{2n+1}_{+}\right)
_{q,\beta}$ will play the roles of the algebras $F_{1}$ and $F_{2}$
(see Section 2) respectively.

(ii) It is easy to get the full list of irreducible $*$-representations
$\pi$ of the algebras ${\cal D}\left(C^{2n+1}\right)_{q,\beta}$ and
${\cal D}\left(H^{2n+1}_{+}\right)_{q,\beta}$ for which the operators
$\pi(x_{1}),...,\pi(x_{n})$ have at least one common eigen-vector.

The principal series of representations is characterized by the condition
$\pi(x_{n}) \not =0$ and are parameterized by the points $(\beta,\varphi)$
of a two-dimensional torus. They are given by (\ref{pi-1})-(\ref{pi-4}).
The representations of degenerate series are either one-dimensional or
obtained from representations for smaller dimensions $n$:
\begin{eqnarray*}
{\bf C}\left[C^{2n+1}\right]_{q}/(z_{n}=z_{n}^{*}=0)
\simeq {\bf C}\left[C^{2n-1}\right]_{q} \\
{\bf C}\left[H^{2n+1}_{+}\right]_{q}/(z_{n}=z_{n}^{*}=0)
\simeq {\bf C}\left[H^{2n-1}_{+}\right]_{q}
\end{eqnarray*}

The limit $q \rightarrow 1$ equips $H^{2n+1}_{+}$ and $C^{2n+1}
\setminus \{0\}$ with a Poisson manifold structure and establishes
a correspondence between the principal series representations and
the symplectic leaves of maximal dimension (${\rm agr} \, z_{n}
=\varphi$). The invariant integral gives the Liouville measure in
the limit, and the degenerate series representations correspond to
the symplectic leaves of dimension $d<2n$ (cf \cite{4}). The parameter
$\beta$, unlike $\varphi$, does not have a simple classical analogue
($q \rightarrow 1$), as $\beta$ distinguishes the homogeneous components
of the quantum $SU(n,1)$-spaces, while in the case $q=1$ the group
$SU(n,1)$ acts transitively on $H^{2n+1}_{+}$ and $C^{2n+1}
\setminus \{0\}$.

In the conclusion of the section, we define the spaces of distributions
and the space of the kernels of integral operators.

Equip the algebras ${\cal D}\left(C^{2n+1}\right)_{q, \beta}$ and
${\cal D}\left(H^{2n+1}_{+}\right)_{q, \beta}$ with the weakest
topology in which are continuous the linear functionals $f \mapsto
f_{IJ}(x_{1},...,x_{n})$ for any $I,J$ and $x_{1},...,x_{n} \in
\mbox{\goth M}_{\beta}$.

The completions ${\cal D}\left(C^{2n+1}\right)_{q, \beta}'$ and
${\cal D}\left(H^{2n+1}_{+}\right)_{q, \beta}'$ will be called
{\em the spaces of distributions\/}. By continuity, they have
a $U_{q}\mbox{\goth sl}(n+1)$-module structures. As follows from
(\ref{1.4}), the pairings
\begin{eqnarray*}
{\cal D}\left(C^{2n+1}\right)_{q, \beta}' \otimes
{\cal D}\left(C^{2n+1}\right)_{q, \beta} \rightarrow {\bf C} \\
{\cal D}\left(H^{2n+1}_{+}\right)_{q, \beta}' \otimes
{\cal D}\left(H^{2n+1}_{+}\right)_{q, \beta} \rightarrow {\bf C} \\
f \otimes \psi \mapsto \int f\psi d\nu \hspace{12mm}
\end{eqnarray*}
yield canonical isomorphisms of the $U_{q}\mbox{\goth sl}(n+1)$-modules
of distributions and the modules dual to the modules of basic functions.

By continuity we define the product of a basic function and a distribution,
which allows to equip ${\cal D}\left(H^{2n+1}_{+}\right)_{q, \beta}'$
with a structure of ${\cal D}\left(H^{2n+1}_{+}\right)_{q, \beta}$-bimodule
and ${\cal D}\left(C^{2n+1}\right)_{q, \beta}'$ with a structure of
${\cal D}\left(C^{2n+1}\right)_{q, \beta}$-bimodule. This structure
is compatible with the structure of $U_{q}\mbox{\goth sl}(n+1)$-module
in the sense of (\ref{1.2}), (\ref{1.3}).

The distributions are decomposable in the series (\ref{2.2}), and they
can be identified with formal series of the form (\ref{2.2}), with the
coefficients $f_{IJ}$ being functions on $\mbox{\goth M}_{\beta}$. The
topology coincides with the topology of the pointwise convergence of
the coefficients $f_{IJ}(x_{1},...,x_{n})$.

For the tensor product ${\bf C}\left[H^{2n+1}_{+}\right]_{q} \otimes
{\bf C}\left[C^{2n+1}\right]_{q}^{\rm op}$ the decomposition (\ref{2.2})
is of the form
\begin{equation}
\label{2.9}
f=\sum_{I \cdot J=0}\sum_{I' \cdot J'=0}
\left(z^{*}\right)^{I}\zeta^{I'}
f_{IJI'J'}(x_{0},...,x_{n};\xi_{0},...,\xi_{n})
z^{J}\left(\zeta^{*}\right)^{J'}
\end{equation}
where $\zeta^{I'}=\zeta_{0}^{i_{0}'}...\zeta_{n}^{i_{n}'}$,
$\left(\zeta^{*}\right)^{J'}=\left(\zeta^{*}_{0}\right)^{j_{0}'}
...\left(\zeta^{*}_{n}\right)^{j_{n}'}$, and
\[ \xi_{0}=\zeta_{0}^{*}\zeta_{0}-\sum_{k>0}\zeta_{k}^{*}\zeta_{k},
\; \; \; \xi_{j}=\sum_{k \geq j}\zeta_{k}^{*}\zeta_{k} \; \; (j \not =0) \]

Passing from the finite sums to the formal series and from the polynomials
$f_{IJI'J'}$ to functions on $\mbox{\goth M}_{\beta} \times \mbox{\goth M}
_{\beta}$, we get the space ${\cal K}\left(H^{2n+1}_{+},C^{2n+1}\right)
_{q,\beta}$ of the kernels of integral operators. In the same way as
above, we equip the space ${\cal K}\left(H^{2n+1}_{+},C^{2n+1}\right)
_{q,\beta}$ with compatible $U_{q}\mbox{\goth sl}(n+1) \otimes
U_{q}\mbox{\goth sl}(n+1)$-module and ${\cal D}\left(H^{2n+1}_{+}
\right)_{q,\beta} \otimes {\cal D}\left(C^{2n+1}\right)_{q,\beta}$-bimodule
structures. We will call {\em intertwining kernels\/} the solutions to
the system (\ref{1.6}) in the space ${\cal K}\left(H^{2n+1}_{+},C^{2n+1}
\right)_{q,\beta}$. As was noted in Section 2, the interwining kernels
correspond to intertwining integral operators ${\cal D}\left(C^{2n+1}
\right)_{q,\beta} \rightarrow {\cal D}\left(H^{2n+1}_{+}\right)_{q,\beta}$.

\section{Non-Integer Powers of the Poisson Kernel}

\setcounter{equation}{0}
\def\theequation{\thesection.\arabic{equation}}

As will be shown, the power $P^{\lambda}$ of the Poisson kernel for
$\lambda \in {\bf Z}_{+}$ can be decomposed into $Q$-binomial ``series'',
with coefficients being rational functions in $q^{2\lambda}$. This
will enable us to perform the analytic continuation in the parameter
$\lambda$. The convergence of the resulting series in the space of
kernels ${\cal K}\left(H^{2n+1}_{+},C^{2n+1}\right)_{q,\beta}$ can
be established in the same way as the following proposition.
\vspace{4mm}

\noindent
{\bf Proposition 4.1.} {\em Suppose that\/} $K'=\sum_{j=1}^{n-1}z_{j}
\zeta_{j}^{*}$ {\em and\/} $K''=\sum_{j=1}^{n-1}q^{-2j}z_{j}^{*}\zeta_{j}$.
{\em Then the series\/}
\begin{equation}
\label{3.1}
\sum_{m=0}^{\infty}\sum_{\begin{array}{c}
m_{1}+m_{2}=m \\ m_{1},m_{2} \geq 0
\end{array}}c_{m_{1}m_{2}}\left(K''\right)^{m_{2}}
\left(K'\right)^{m_{1}}
\end{equation}
{\em converges in\/} ${\cal K}\left(H^{2n+1}_{+},C^{2n+1}\right)_{q,\beta}$
{\em for any\/} $c_{m_{1}m_{2}} \in {\bf C}$.
\vspace{4mm}

\noindent
{\em Proof\/}. Choosing a point $(x_{1},...,x_{n};\xi_{1},...,\xi_{n})
\in \mbox{\goth M}_{\beta} \times \mbox{\goth M}_{\beta}$ and a quadruple
of multi-indices $(I,J,I',J')$, we bring the terms of the series (\ref{3.1})
to the ``normal form'' (\ref{2.9}) by using the commutation relations
(\ref{1.9}),(\ref{1.10}). Then it is suffices to show that for some $M \in
{\bf Z}_{+}$ the terms with the indices $m>M$ give the zero contribution
to $f_{IJI'J'}(x_{0},...,x_{n};\xi_{0},...,\xi_{n})$. But this follows
from the definition of $\mbox{\goth M}_{\beta}$ and the following
identities in the algebra ${\bf C}_{q,q^{-1}}^{n+1} \otimes \left(
{\bf C}_{q,q^{-1}}^{n+1}\right)^{\rm op}$:
\begin{eqnarray*}
\left(z_{1}^{*}\right)^{m_{1}}\left(z_{2}^{*}\right)^{m_{2}}
...\left(z_{n-1}^{*}\right)^{m_{n-1}}z_{1}^{m_{1}}z_{2}^{m_{2}}
...z_{n-1}^{m_{n-1}}= \hspace{1mm} \\
=const_{1}(m_{1},...,m_{n-1})\prod_{k=1}^{n-1}\prod_{j=1}^{m_{k}}
\left(x_{k}-q^{-2}x_{k+1}\right), \\
\zeta_{1}^{m_{1}}\zeta_{2}^{m_{2}}...\zeta_{n-1}^{m_{n-1}}
\left(\zeta_{1}^{*}\right)^{m_{1}}\left(\zeta_{2}^{*}\right)
^{m_{2}}...\left(\zeta_{n-1}^{*}\right)^{m_{n-1}}= \hspace{1mm} \\
=const_{2}(m_{1},...,m_{n-1})\prod_{k=1}^{n-1}\prod_{j=1}^{m_{k}}
\left(\xi_{k}-q^{-2}\xi_{k+1}\right) \hspace{2mm}
\end{eqnarray*}

By (\ref{1.18}), for $\lambda \in {\bf Z}_{+}$ the powers of
$P=K_{1}K_{2}$ are given by
\begin{equation}
\label{3.2}
P^{\lambda}=q^{\lambda(\lambda+1)}\left(z_{0}^{*}\zeta_{0}-K''-
q^{-2n}z_{n}^{*}\zeta_{n}\right)^{\lambda}\left(z_{0}\zeta_{0}^{*}
-K'-z_{n}\zeta_{n}^{*}\right)^{\lambda}
\end{equation}

Note that the summands in any of the parentheses quasi-commute.
For instance, we have that
\begin{eqnarray*}
\left(z_{0}\zeta_{0}^{*}\right)K'=q^{2}K'\left(z_{0}\zeta_{0}^{*}\right), &
K'\left(z_{n}\zeta_{n}^{*}\right)=q^{2}\left(z_{n}\zeta_{n}^{*}\right)K' \\
\left(z_{0}^{*}\zeta_{0}\right)K''=q^{-2}K''\left(z_{0}^{*}\zeta_{0}\right), &
K''\left(z_{n}^{*}\zeta_{n}\right)=q^{-2}\left(z_{n}^{*}\zeta_{n}\right)K''
\end{eqnarray*}

This allows us to use the well-known $q$-analogue
of the Newton binomial formula:
\vspace{4mm}

\noindent
{\bf Proposition 4.2.} {\em If\/} $a$ {\em and\/} $b$ {\em are elements
of an associative algebra such that\/} $ab=q^{2}ba$, {\em then\/}
\begin{eqnarray*}
(a+b)^{m}=\sum_{j=0}^{m}\frac{\left(q^{2};q^{2}\right)_{m}}
{\left(q^{2};q^{2}\right)_{j}\left(q^{2};q^{2}\right)_{m-j}}
b^{j}a^{m-j}, \\
\left(t;q^{2}\right)_{k} \stackrel{\rm def}{=}
(1-t)(1-q^{2}t)...(1-q^{2(k-1)}t) \hspace{2mm}
\end{eqnarray*}

\noindent
{\bf Corollary 4.3.} {\em We have that\/}
\begin{eqnarray*}
\left(z_{0}\zeta_{0}^{*}-K'-z_{n}\zeta_{n}^{*}\right)^{\lambda}
  & = & \sum_{k}\frac{\left(q^{2\lambda};q^{-2}\right)_{k}}
{\left(q^{2};q^{2}\right)_{k}} \sum_{i}\frac{\left(q^{2(\lambda
+k)};q^{-2}\right)_{i}}{\left(q^{2};q^{2}\right)_{i}} \cdot \\
  & \cdot & q^{-2(i+k)(\lambda-i-k)}\left(z_{0}\zeta_{0}^{*}
\right)^{\lambda-i-k}\left(-z_{n}\zeta_{n}^{*}\right)^{i}
\left(-K'\right)^{k}, \\
\left(z_{0}^{*}\zeta_{0}-K''-q^{-2n}z_{n}^{*}\zeta_{n}\right)^{\lambda}
  & = & \sum_{k}\frac{\left(q^{-2\lambda};q^{2}\right)_{k}}
{\left(q^{-2};q^{-2}\right)_{k}} \sum_{i}\frac{\left(q^{-2(\lambda-k)};
q^{2}\right)_{i}}{\left(q^{-2};q^{-2}\right)_{i}} \cdot \\
  & \cdot & q^{2ki}\left(-K''\right)^{k}\left(-q^{-2n}z_{n}^{*}
\zeta_{n}\right)^{i}\left(z_{0}^{*}\zeta_{0}\right)^{\lambda-i-k}
\end{eqnarray*}

It follows from (\ref{3.2}) and Corollary 4.3 that $P^{\lambda}$
can be decomposed into a sum
\begin{equation}
\label{3.3}
P^{\lambda}=\sum_{l_{1},l_{2}} \left(-K''\right)^{l_{1}}
f_{l_{1}l_{2}}\left(q^{2\lambda}\right) \left(K'\right)^{l_{2}}
\end{equation}
of the form (\ref{3.1}), where
\begin{eqnarray}
\label{3.4}
f_{l_{1}l_{2}}\left(q^{2\lambda}\right) & = &
\sum_{i_{1},i_{2}}q^{\lambda(\lambda+1)}
const(\lambda,i_{1},i_{2},l_{1},l_{2}) \cdot \\
\nonumber
  & \cdot & \left(-q^{-2n}z_{n}^{*}\zeta_{n}\right)^{i_{1}}
\left(z_{0}^{*}\zeta_{0}\right)^{\lambda-l_{1}-i_{1}}
\left(z_{0}\zeta_{0}^{*}\right)^{\lambda-l_{2}-i_{2}}
\left(-z_{n}\zeta_{n}^{*}\right)^{i_{2}} \\
\nonumber
const(\lambda,i_{1},i_{2},l_{1},l_{2}) & = &
q^{l_{1}(l_{1}+1)-l_{2}(l_{2}+1)}\frac{\left(q^{-2\lambda};
q^{2}\right)_{l_{1}}\left(q^{-2\lambda};q^{2}\right)_{l_{2}}}
{\left(q^{2};q^{2}\right)_{l_{1}}\left(q^{2};q^{2}\right)_{l_{2}}} \cdot \\
\nonumber
  & \cdot & (-1)^{l_{1}+l_{2}+i_{1}+i_{2}}q^{i_{1}(i_{1}+1)-
i_{2}(i_{2}+1)+2i_{1}l_{1}+2i_{2}l_{2}+2(i_{2}+l_{2})^{2}} \cdot \\
\label{3.5}
  & \cdot & \frac{\left(q^{-2(\lambda-l_{1})};q^{2}\right)_{i_{1}}
\left(q^{-2(\lambda+l_{1})};q^{2}\right)_{i_{2}}}
{\left(q^{2};q^{2}\right)_{i_{1}}\left(q^{2};q^{2}\right)_{i_{2}}}
\end{eqnarray}

\noindent
{\bf Proposition 4.4.} {\em The powers of the Poisson kernel\/}
$P^{\lambda}$ {\em have a decomposition of the form\/} (\ref{2.9}),
{\em with the coefficients of the power series\/} $\xi_{1}^{-\lambda}
f_{IJI'J'}$ {\em being polynomials in\/} $q^{2\lambda}$, $q^{-2\lambda}$
({\em in particular, rational functions in\/} $q^{2\lambda}$).
\vspace{4mm}

To prove it, it suffices to go to a ``normal'' ordering of the generators
in $P^{\lambda}$ (as in (\ref{2.9})) by using (\ref{3.3})-(\ref{3.5})
and the following simple statement.
\vspace{4mm}

\noindent
{\bf Proposition 4.5.} {\em Suppose that\/} $m_{1} \leq \lambda$,
$m_{2} \leq \lambda$, $m={\rm max}(m_{1},m_{2})$. {\em Then\/}
\begin{eqnarray}
\label{3.6}
q^{\lambda(\lambda+1)}\left(z_{0}^{*}\zeta_{0}\right)^{\lambda-m_{1}}
\left(z_{0}\zeta_{0}^{*}\right)^{\lambda-m_{2}} & = &
q^{m(2\lambda-m+1)}\left(z_{0}^{*}\zeta_{0}\right)^{m-m_{1}} \cdot \\
\nonumber
  & \cdot & \frac{\left(-q^{-2(\lambda-m)}x_{1};q^{2}\right)_{\infty}}
{\left(-x_{1};q^{2}\right)_{\infty}}\xi_{1}^{\lambda-m}
\left(z_{0}\zeta_{0}^{*}\right)^{m-m_{2}}
\end{eqnarray}
{\em where\/}
\[ \left(t;q^{2}\right)_{\infty}=\prod_{j=0}^{\infty}
\left(1-q^{2j}t\right)=\sum_{l=1}^{\infty}
\frac{t^{l}q^{l(l-1)}}{\left(q^{2};q^{2}\right)_{l}} \]

\noindent
{\em Remark\/}. The crucial point in the proof of Proposition 4.5
is the elimination of the multiple $q^{\lambda^{2}}$ due to the
fact that $\zeta_{0}^{\lambda-m}\left(\zeta_{0}^{*}\right)^{\lambda
-m}=q^{-(\lambda+m)(\lambda+m+1)}\xi_{1}^{\lambda-m}$.
\vspace{4mm}

Now, by analytical continuation, we can define $P^{\lambda}$ as
a power series for any $\lambda \in {\bf R}$. The formulas
(\ref{3.3})-(\ref{3.6}) hold in this case. Let us bring the
coefficients $f_{l_{1}l_{2}}$ of the series (\ref{3.3}) to the
form (\ref{2.9}):
\begin{eqnarray}
\label{3.7}
f_{l_{1}l_{2}} & = & \sum_{\begin{array}{c}
j_{0},j_{n},k_{0},k_{n} \\ j_{0}k_{0}=j_{n}k_{n}=0
\end{array}} \left(z_{0}^{*}\zeta_{0}\right)^{j_{0}}
\left(-q^{-2n}z_{n}^{*}\zeta_{n}\right)^{j_{n}}
\xi_{1}^{\lambda} \cdot \\
\nonumber
  & \cdot & f_{l_{1}l_{2}j_{0}j_{n}k_{0}k_{n}}
\left(q^{2\lambda};x_{1},x_{n},\xi_{1}^{-1}\xi_{n}\right)
\left(z_{0}\zeta_{0}^{*}\right)^{k_{0}}
\left(z_{n}\zeta_{n}^{*}\right)^{k_{n}}
\end{eqnarray}

As can be seen from the proof of Proposition 4.1, in
order to get a convergent series in the space of kernels
${\cal K}\left(H^{2n+1}_{+},C^{2n+1}\right)_{q,\beta}$ it
suffices to establish convergence of any of the power
series $f_{l_{1}l_{2}j_{0}j_{n}k_{0}k_{n}}$ in some
neighbourhood of zero and to continue them analytically.

Recall the definition of the {\em basic hypergeometric series\/}
$_{r+1}\Phi_{r+j}$ (cf. \cite{5}):
\begin{eqnarray*}
_{r+1}\Phi_{r+j}\left(\begin{array}{c}
a_{0},...,a_{r} \\ b_{1},...,b_{r+j}
\end{array};q,x\right)\stackrel{\rm def}{=} \hspace{10mm} \\
\stackrel{\rm def}{=}\sum_{m=0}^{\infty}
\frac{(a_{0};q)_{m}...(a_{r};q)_{m}(-1)^{jm}q^{ij(m-1)/2}}
{(b_{1};q)_{m}...(b_{r+j};q)_{m}(q;q)_{m}}x^{m}
\end{eqnarray*}

\noindent
{\bf Proposition 4.6.} {\em we have that\/}
\begin{eqnarray}
\label{3.8}
f_{000000}\left(q^{2\lambda};x_{1},x_{n},
\xi_{1}^{-1}\xi_{n}\right)= \hspace{30mm} \\
\nonumber
=\frac{\left(-q^{-2\lambda}x_{1};q^{2}\right)_{\infty}}
{\left(-x_{1};q^{2}\right)_{\infty}} \, _{2}\Phi_{2}
\left(\begin{array}{cc}
q^{-2\lambda}, & q^{-2\lambda} \\
q^{2}, & -q^{-2\lambda}x_{1}
\end{array};q^{2},-q^{2(\lambda+n-2)}
x_{n}\xi_{1}^{-1}\xi_{n}\right)
\end{eqnarray}

The proof is straightforward and based on (\ref{3.4}),
(\ref{3.5}) and Corollary 4.5.
\vspace{4mm}

\noindent
{\bf Proposition 4.7.} {\em For any\/} $\lambda \in {\bf R}$,
{\em the series\/} $f_{l_{1}l_{2}j_{0}j_{n}k_{0}k_{n}}$
{\em converges in some neighbourhood of zero and can be
continued analytically in the region\/} $x_{1}>0$.
\vspace{4mm}

In the special case $l_{1}=l_{2}=j_{0}=j_{n}=k_{0}=k_{n}=0$
this proposition follows from the previous one. In the general
case it can de proved in the way analogous to that of the
proof of Proposition 4.6. Namely, the analytic continuation
can be achieved by applying (\ref{3.6}), while the convergence
is provided by the multiple $q^{i_{1}^{2}+i_{2}^{2}}$ in (\ref{3.5}).
\vspace{4mm}

\noindent
{\bf Definition 4.8.} For any real $\lambda$ the formulas (\ref{3.3}),
(\ref{3.7}) and Proposition 4.7 determine a convergent series in the
space of kernels ${\cal K}\left(H^{2n+1}_{+},C^{2n+1}\right)_{q,\beta}$.
We call its sum {\em the\/} $\lambda$-{\em th power of the Poisson
kernel\/} and denote it by $P^{\lambda}$.
\vspace{4mm}

Following (\ref{1.8}), we define an anti-linear anti-involution
in the space of kernels by
\begin{eqnarray*}
z_{j} \mapsto z_{j}^{*}, \; \; \;
f(x_{1},...,x_{n};\xi_{1},...,\xi_{n})
\mapsto f(x_{1},...,x_{n};\xi_{1},...,\xi_{n}), \\
\zeta_{j} \mapsto (\xi_{1}...\xi_{n})^{-1}\zeta_{j}^{*}
(\xi_{1}...\xi_{n})=q^{2(j-n)}\xi_{j}^{*} \hspace{16mm}
\end{eqnarray*}

It is clear that $\left(K''\right)^{*}=q^{-2n}K'$. Therefore,
$P^{\lambda}=\left(P^{\lambda}\right)^{*}$.
\vspace{4mm}

\noindent
{\bf Proposition 4.9.} {\em For any\/} $\lambda \in {\bf R}$,
{\em the kernal\/} $P^{\lambda} \in {\cal K}\left(H^{2n+1}_{+},
C^{2n+1}\right)_{q,\beta}$ {\em is intertwining\/}.
\vspace{4mm}

\noindent
{\em Proof\/}. Consider the kernel $L(\lambda)=\left(\xi \otimes
1-1 \otimes S^{-1}(\xi)\right)\left(P^{\lambda}\right)$. By
(\ref{1.9}), (\ref{1.10}), (\ref{2.3})-(\ref{2.6}),
(\ref{1.12})-(\ref{1.17}), the coefficients in the decomposition
of $L(\lambda)$ in a series of the form (\ref{2.9}) are
analytical functions in $x_{1},...,x_{n},\xi_{1}^{-1}\xi_{2},
...,\xi_{1}^{-1}\xi_{n}$ in the region $x_{1}>0$. Close to zero
they can be decomposed into power series, with coefficients
being rational functions in $q^{2\lambda}$. But at the points
$q^{2},q^{4},q^{6},...$ these rational functions have zeroes,
since the kernel $P$ and hence the kernels $P^{\lambda}$
($\lambda \in {\bf Z}_{+}$) are intertwining. Therefore,
$L(\lambda) \equiv 0$. The proposition is proved.
\vspace{4mm}

Analogously, one can prove that $P^{\lambda}P^{m}=P^{m}P^{\lambda}=
P^{\lambda+m}$ for any $\lambda \in {\bf R}$, $m \in {\bf Z}_{+}$.

Some applications of the intertwining kernals are related to the
fact that they are generating functions for spherical functions
on the quantum homogeneous spaces.
\vspace{4mm}

\noindent
{\em Example\/}. When $n=1$, $\lambda=-l-1$, $l \in {\bf Z}_{+}$,
the function $f_{000000}$ is a zonal spherical function
corresponding to the spin $l$ finite-dinemsional representation
of $U_{q}\mbox{\goth sl}(2)$. By (\ref{3.8}), we get that
\begin{eqnarray}
\nonumber
f_{000000} & = & \frac{\left(-q^{2(l+1)}x_{1};q^{2}\right)_{\infty}}
{\left(-x_{1};q^{2}\right)_{\infty}} \, _{2}\Phi_{2}
\left(\begin{array}{cc}
q^{2(l+1)}, & q^{2(l+1)} \\
q^{2}, & q^{2(l+1)}x_{1}
\end{array};q^{2},-q^{-2l}x_{1}\right)= \\
\label{3.9}
  & = & \, _{2}\Phi_{1}
\left(\begin{array}{c}
q^{2(l+1)}, q^{-2l} \\ q^{2}
\end{array};q^{2},-qx\right)
\end{eqnarray}

Here we have used the $q$-analogue of the Pfaff transform
(cf. (1.32) in \cite{5}). The right hand side of (\ref{3.9})
is the zonal spherical function (cf. \cite{6}). The generating
function for the $q$-analogues of the Clebsch-Gordan coefficients
obtained in \cite{7} also is an intertwining kernel in the
sense of Section 2 of the present paper.

\section{Quantum Radon Transform}

\setcounter{equation}{0}
\def\theequation{\thesection.\arabic{equation}}

We are interested in functions of the Poisson kernel which are
not necessarily polynomials. In Section 4 were defined the powers
$P^{\lambda} \in {\cal K}\left(H^{2n+1}_{+},C^{2n+1}\right)_{q,\beta}$
for any $\lambda \in {\bf R}$. In the present section we will obtain
their decomposition of the form
\begin{equation}
\label{4.1}
P^{\lambda}=\sum_{j \in {\bf Z}}q^{2(j+\beta)\lambda}R_{j}
\end{equation}
where $R_{j} \in {\cal K}\left(H^{2n+1}_{+},C^{2n+1}\right)_{q,\beta}$.
The formula (\ref{4.1}) allows to interpret the kernels $R_{j}$ as
$\delta$-functions $\delta\left(P-q^{2(j+\beta)}\right)$ and the
integral transforms with these kernels as quantum analogues of the
Radon transform.

In Section 4 we were interested in the functions $f_{l_{1}l_{2}j_{0}
j_{n}k_{0}k_{n}}$ for a fixed $\lambda$. Now we fix $(x_{1},...,x_{n})
\in \mbox{\goth M}_{\beta}$, $(\xi_{1},...,\xi_{n}) \in \mbox{\goth M}
_{\beta}$.
\vspace{4mm}

\noindent
{\bf Proposition 5.1.} {\em After the change of the variable\/}
$q^{2\lambda}=u$ {\em the function\/} $f_{l_{1}l_{2}j_{0}j_{n}k_{0}
k_{n}}(u)$ {\em is holomorphic in some neighbourhood of zero but
the point\/} $u=0$ {\em for any\/} $l_{1},l_{2},j_{0},j_{n},k_{0},
k_{n}$.
\vspace{4mm}

This statement follows from (\ref{3.5}), (\ref{3.6}) and the
commutation relations between the generators of the algebra
${\bf C}_{q,q^{-1}}^{n+1} \otimes \left({\bf C}_{q,q^{-1}}^{n+1}
\right)^{\rm op}$.
\vspace{4mm}

\noindent
{\bf Corollary 5.2.} {\em There exist such functions\/} $r_{l_{1}l_{2}
j_{0}j_{n}k_{0}k_{n}}^{(m)}$ {\em on\/} $\mbox{\goth M}_{\beta} \times
\mbox{\goth M}_{\beta}$, {\em that\/}
\begin{equation}
\label{4.2}
f_{l_{1}l_{2}j_{0}j_{n}k_{0}k_{n}}=\sum_{m \in {\bf Z}}
r^{(m)}_{l_{1}l_{2}j_{0}j_{n}k_{0}k_{n}}q^{2m\lambda}
\end{equation}

Plugging (\ref{4.2}) into (\ref{3.7}), we get that the dependence on
$\lambda$ is reflected in the factor $\left(\xi_{1}q^{2m}\right)
^{\lambda}$ only.  Therefore, only those values of $\xi_{1}$ make
a contribution into $\delta(P-a)$ for which $\xi_{1}=aq^{-2m}$.
\vspace{4mm}

\noindent
{\bf Definition 5.3.} {\em The kernel\/} $\delta\left(P-q^{2\beta}
\right)$ {\em of the quantum Radon transform\/} is the sum of the
series (\ref{3.3}) where
\begin{eqnarray*}
f_{l_{1}l_{2}} & = & \sum_{j_{0}k_{0}=j_{n}k_{n}=0}\left(z_{0}^{*}
\zeta_{0}\right)^{j_{0}}\left(-q^{-2n}z_{n}^{*}\zeta_{n}\right)
^{j_{n}} \cdot \\
  & \cdot & r^{(m)}_{l_{1}l_{2}j_{0}j_{n}k_{0}k_{n}}|_{m=
\frac{1}{2}\log_{q}\xi-\beta}\left(z_{0}\zeta_{0}^{*}\right)
^{k_{0}}\left(z_{n}\zeta_{n}^{*}\right)^{k_{n}}
\end{eqnarray*}

The convergence of the series (\ref{3.3}) in the space ${\cal K}
\left(H^{2n+1}_{+},C^{2n+1}\right)_{q,\beta}$ can be proved in
the same way as Proposition 4.1. Similarly to the above definition,
one can introduce the kernels $R_{j}=\delta\left(P-q^{2(j+\beta)}
\right)$ for any integer $j$. The decomposition (\ref{4.1}) follows
from the above constructions.
\vspace{4mm}

\noindent
{\bf Proposition 5.4.} {\em The kernel\/} $\delta\left(P-q^{2\beta}
\right)$ {\em is intertwining\/}.
\vspace{4mm}

This follows from the formula
\[ r^{(m)}_{l_{1}l_{2}j_{0}j_{n}k_{0}k_{n}}=\frac{1}{2\pi i}
\int_{\begin{picture}(1,10)(0,0)
\put(0.5,-1){\line(0,1){6}}
\end{picture} \, u \,
\begin{picture}(1,10)(0,0)
\put(0.5,-1){\line(0,1){6}}
\end{picture}=\varepsilon}
\frac{1}{u^{m+1}}f_{l_{1}l_{2}j_{0}j_{n}k_{0}k_{n}}du \]
and from the fact that the action of $U_{q}\mbox{\goth sl}(n+1)
\otimes U_{1}\mbox{\goth sl}(n+1)$ commutes with the integration.
Indeed, repeating for small
$\begin{picture}(1,10)(0,0)
\put(0.5,-2){\line(0,1){8}}
\end{picture} \, u
\begin{picture}(1,10)(0,0)
\put(0.5,-2){\line(0,1){8}}
\end{picture} \, =q^{\mbox{\rm \footnotesize Re} \, \lambda}$
the prrof of Proposition 4.9, we get that the functions under the
integration symbol -- which are obtained by applying the operator
$\xi \otimes 1-1 \otimes S^{-1}(\xi)$ to the kernel $\delta\left(
P-q^{2\beta}\right)$ -- are equal to zero.
\vspace{4mm}

\noindent
{\em Remark\/}. Suppose that $m \in {\bf Z}_{+}$. Then $P^{\lambda}
P^{m}=P^{m}P^{\lambda}=P^{\lambda+m}$ implies that $P^{m}\delta
\left(P-q^{2\beta}\right)=\delta\left(P-q^{2\beta}\right)P^{m}=
q^{2\beta m}P^{m}$.

\section{Intertwining Kernels and Hypergeometric
        \newline Functions}

\setcounter{equation}{0}
\def\theequation{\thesection.\arabic{equation}}

Recall that the algebra of intertwining kernels is a subalgebra
of $F_{2} \otimes F_{1}^{\rm op}$ and is defined by (\ref{1.6}).
In what follows the role of $F_{1}$ is played by ${\bf C}\left[
C^{2n+1}\right]_{q}$ (cf. Section 3), and the role of $F_{2}$ by
the algebra $\left(\oplus_{i=1}^{N} {\bf C}^{n+1}\right)_{q,q^{-1}}
^{\rm op}$ defined below.

We will use the $R$-matrix notation of \cite{3}. Let $e_{ij} \in
\mbox{\rm Mat}_{n+1}({\bf C})$ be the matrix unit and
\begin{eqnarray*}
R' & = & q^{2}\sum_{i}e_{ii} \otimes e_{ii}+q\sum_{i \not =j}
e_{ij} \otimes e_{ji}+\left(q^{2}-1\right)\sum_{i<j}e_{ii}
\otimes e_{jj}, \\
R'' & = & q^{-1}\sum_{i}e_{ii} \otimes e_{ii}+q\sum_{i \not =j}
e_{ij} \otimes e_{ji}+\left(q-q^{-1}\right)\sum_{i<j}q^{i-j}e_{ij}
\otimes e_{ij}.
\end{eqnarray*}

Define $\left(\oplus_{i=1}^{N} {\bf C}^{n+1}\right)_{q,q^{-1}}
^{\rm op}$ as the $*$-algebra generated by $z_{ij}$ ($1 \leq i
\leq N$, $0 \leq j \leq n$) and the relations
\begin{eqnarray*}
\sum_{i_{1}j_{1}}R_{iji_{1}j_{1}}'z_{bi_{1}}z_{aj_{1}}=
qz_{ai}z_{bj} \; \; (a<b), \hspace{3mm} \\
\sum_{i_{1}j_{1}}R_{iji_{1}j_{1}}''\left(\varepsilon_{i_{1}}
z_{bi_{1}}^{*}\right)z_{aj_{1}}=q^{-1}z_{ai}\left(
\varepsilon_{j}z_{bj}^{*}\right),
\end{eqnarray*}
where $\varepsilon_{i}=(-1)^{\delta_{i0}}q^{-i}$.

Define an action of the generators $X_{i}^{\pm}$, $K_{i}^{\pm}$
of $U_{q}\mbox{\goth su}(n,1)$ on the generators $z_{ij}$ of
$\left(\oplus_{i=1}^{N} {\bf C}^{n+1}\right)_{q,q^{-1}}^{\rm op}$
by (\ref{1.12})-(\ref{1.14}), replacing $z_{j},z_{j-1},z_{j+1}$
by $z_{ij},z_{i,j-1},z_{i,j+1}$ respectively. Extend this action
so that $\left(\oplus_{i=1}^{N} {\bf C}^{n+1}\right)_{q,q^{-1}}
^{\rm op}$ becomes a $U_{q}\mbox{\goth su}(n,1)$-module $*$-algebra.

The uniqueness of such an extension is obvious due to the conditions
(\ref{1.2}), (\ref{1.3}), (\ref{1.7}). The existence follows from
the fact that the matrices $R'$ and $R''$ correspond to linear
operators which intertwine certain representations of $U_{q}
\mbox{\goth sl}(n+1)$. (For any $a$ the elements $z_{aj}$ form
a standard basis of the space of the vector representation, while
the elements $\varepsilon_{j}z_{aj}^{*}$ form a standard basis of
the space of the covector one (cf. \cite{2}).)

Note that for $N=1$ the resulting $U_{q}\mbox{\goth su}(n,1)$-module
$*$-algebra is isomorphic to ${\bf C}_{q,q^{-1}}^{n+1}$.

We will identify the algebras $F_{1}^{\rm op}={\bf C}\left[C^{2n+1}
\right]_{q}^{\rm op}$ and $F_{2}=\left(\oplus_{i=1}^{N}{\bf C}^{n+1}
\right)_{q,q^{-1}}$ with their images under the canonical embeddings
into $F_{2} \otimes F_{1}^{\rm op}$. Following (\ref{1.7}), (\ref{2.8}),
introduce an anti-linear involution in $F_{2} \otimes F_{1}^{\rm op}$ by
\[ z_{ij} \mapsto z_{ij}^{*}, \; \; \;
\zeta_{j} \mapsto q^{2(j-n)}\zeta_{j}^{*} \]

The following proposition is a straightforward consequence
of the definitions.
\vspace{4mm}

\noindent
{\bf Proposition 6.1.} {\em The elements\/}
\[ K_{i}=z_{i0}\zeta_{0}^{*}-\sum_{j=1}^{n}z_{ij}\zeta_{j}^{*} \]
{\em of\/} $F_{2} \otimes F_{1}^{\rm op}$ {\em are intertwining
kernels, and we have that\/}
\begin{eqnarray*}
K_{i}K_{j}=qK_{j}K_{i} \; \; (i<j), \; \; \;
K_{i}K_{j}^{*}=qK_{j}^{*}K_{i} \; \; (i \not =j), \\
K_{i}K_{i}^{*}=q^{2}K_{i}^{*}K_{i} \hspace{27mm}
\end{eqnarray*}

\noindent
{\bf Corollsry 6.2.} {\em The Poisson kernels\/}
\[ P_{i}=q^{-2n}K_{i}K_{i}^{*} \]
{\em are intertwining for any\/} $1 \leq i \leq n$.
{\em Moreover, we have that for any\/} $i,j$,
\[ P_{i}P_{j}=P_{j}P_{i}, \; \; \; P_{i}^{*}=P_{i} \]

As was explained in Section 4, the powers of the Poisson
kernel are generating functions for some polynomials of
hypergeometric form. This allows to consider the intertwining
kernels $\prod_{j=1}^{N}P_{j}^{l_{j}}$ ($l_{j} \in {\bf Z}_{+}$)
as generating functions for elements of the algebra $F_{2}$
which generalize classical orthogonal polynomials.

Repeating the considerations of Section 4, one can get rid of
the restriction $l_{1} \in {\bf Z}_{+}$ for one of the powers
of the Poisson kernels. This enables us to incorporate the case
$\sum_{j=1}^{n}l_{j}=-n$ and obtain some invariants in $\left(
\oplus_{i=1}^{N}{\bf C}^{n+1}\right)_{q,q^{-1}}$ by integrating
the kernel $\prod_{j=1}^{N}P_{j}^{l_{j}}$.

This approach considered in the limit case $q=1$ is analogous
to the approach of V.A.Vasiliev, I.M.Gelfand, and A.B.Zelevinsky
to generalized hypergeometric functions (cf. \cite{8}).
\vspace{5mm}

\Large

\noindent
{\bf Acknowledgments.}
\vspace{3.5 mm}

\normalsize

\noindent
I want to express my gratitude to V.G.Drinfeld for fruitful
discussions of the role of intertwining operators in
representation theory and harmonic analysis.

\end{document}